\input form
%

\def\lie#1{L_{#1}}
\def\liec#1{L_{#1}}

\def\diag{\mathop {\rm diag}}
\def\eqalignno#1{\leqalignno{#1}}

\def\ordnorma#1{\|{#1}\|}
\def\bignorma#1{\bigl\|{#1}\bigr\|}

\def\biggnorma#1{\biggl\|{#1}\biggr\|}
\def\cbin#1#2{{{#1}\choose{#2}}}

\def\uno{{\mat 1}}
\def\ilt{\mathrel{\triangleleft}}

\cita{Arnold-1980}{V.I.\ Arnold: {\it Geometrical Methods in the Theory of
Ordinary Differential Equations}, Springer Verlag, (1980)}

\cita{Buff-2006}{Buff, X. and Ch\'eritat, A.: {\it The Brjuno
function continuously estimates the size of quadratic Siegel disks},
Ann. of Math. {\bf 164}, 265--312 (2006).} 

\cita{Carletti-Marmi-2000}{Carletti, T. and Marmi, S.: {\it
Linearization of analytic and non-analytic germs of diffeomorphisms
of $(\complessi, 0)$}, Bull. Soc. Math. France, {\bf 128}, 69--85
(2000).}

\cita{Giorgilli-1997.4}{Giorgilli, A. and Locatelli, U.: {\it On
classical series expansions for quasi-periodic motions}, MPEJ {\bf 3}
N. 5 (1997).}

\cita{Giorgilli-1999}{A. Giorgilli and U. Locatelli: {\it A classical
self-contained proof of Kolmogorov's theorem on invariant tori}, in
{\it Hamiltonian systems with three or more degrees of freedom},
Carles Sim\'o ed., NATO ASI series C, Vol. 533, Kluwer Academic
Publishers, Dordrecht--Boston--London (1999).}

\cita{Giorgilli-2001}{Giorgilli, A.: {\it Unstable equilibria of
Hamiltonian systems}, Disc. and Cont. Dynamical Systems, Vol. 7, N. 4,
855--871 (2001).}

\cita{Giorgilli-2003.1}{A.\ Giorgilli: {\it Notes on exponential stability
of Hamiltonian systems}, in {\it Dynamical Systems, Part I:
Hamiltonian systems and Celestial Mechanics} Pubblicazioni del Centro
di Ricerca Matematica Ennio De Giorgi, Pisa, 87--198 (2003).}

\cita{Gio-Mur-2006} {A.\ Giorgilli and D.\ Muraro: {\it Exponentially
stable manifolds in the neighbourhood of elliptic equilibria},
Boll. Unione Mat. Ital. Sez. B {\bf 9}, 1--20 (2006).}

\cita{Gio-Mar-2010}{A.\ Giorgilli and S.\ Marmi: {\sl Convergence radius
in the Poincar\'e-Siegel problem}, DCDS Series S, {\bf 3}, 601--621
(2010).}

\cita{Giorgilli-2012}{A.\ Giorgilli: {\it On a theorem of Lyapounov},
Rendiconti dell'Istituto Lombardo Accademia di Scienze e Lettere,
Classe di Scienze Matematiche e Naturali, to appear (see also
arXiv:1303.7322).}

\cita{Giorgilli-2013}{A.\ Giorgilli: {\it On the representation of
maps by Lie transforms}, Rendiconti dell'Istituto Lombardo Accademia
di Scienze e Lettere, Classe di Scienze Matematiche e Naturali, to
appear (see also arXiv:1211.5674).} 

\cita{Gio-Loc-San-2013}{A.\ Giorgilli, U.\ Locatelli and M.\
Sansottera: {\it On the convergence of an algorithm constructing the
normal form for lower dimensional elliptic tori in planetary systems},
(arXiv:1401.6529).}

\cita{Groebner-1960}{W.\ Gr\"obner: {\it Die  Lie-Reihen und Ihre
Anwendungen}, VEB Deutscher Verlag der Wissenschaften, Mathematische 
Monographien, {\bf 3} (1960).
Italian translation: {\it Serie di Lie e loro applicazioni},
Ed. Cremonese, Roma (1973).}

\cita{Schroeder-1871}{E. Schr\"oder: {\it \"Uber iterierte
Functionen}, Math. Ann. {\bf 3}, 296-322 (1871).}

\cita{Siegel-1942}{Siegel, C.L.: {\it Iterations of analytic functions},
Annals of Math. {\bf 43}, 607--612 (1942).}

\cita{Yoccoz-1995}{Yoccoz, J.-C.: {\it Th\'eor\`eme de Siegel, nombres de Bruno et
polyn\^omes quadratiques}  Ast\'erisque {\bf 231} 3--88 (1995).}

\cita{Yoccoz-2002}{Yoccoz, J.-C.: {\it Analytic linearisation of circle
diffeomorphisms} in ``Dynamical Systems and Small Divisors'',
Lecture Notes in Mathematics {\bf 1784} 125--173 (2002).}

\def\testos{A.\ Giorgilli, U.\ Locatelli and M.\ Sansottera}
\def\testod{Improved convergence estimates for the Schr\"oder-Siegel problem}

\title{Improved convergence estimates\riga for the Schr\"oder-Siegel problem}

\author{\it ANTONIO GIORGILLI
\hfill\break Dipartimento di Matematica,
Universit\`a degli Studi di Milano,\hfill\break
Via Saldini 50, 20133\ ---\  Milano, Italy.}

\author{\it UGO LOCATELLI\hfil\break
Dipartimento di Matematica, 
Universit\`a degli Studi di Roma ``Tor Vergata'',\hfil\break 
Via della Ricerca Scientifica 1, 00133\ ---\ Roma, Italy.}

\author{\it MARCO SANSOTTERA\hfil\break 
Universit\`a degli Studi di Milano,\hfill\break
Via Saldini 50, 20133\ ---\  Milano, Italy,\hfil\break
D\'epartement de Math\'ematique, University of Namur \& NAXYS,\hfil\break
Rempart de la Vierge 8, B-5000\ ---\ Namur, Belgium.}

\abstract{We reconsider the Schr\"oder-Siegel problem
of conjugating an analytic map in $\complessi$ in the neighborhood of
a fixed point to its linear part, extending it to the case of
dimension $n\gt 1$.  Assuming a condition which is equivalent to
Bruno's one on the eigenvalues $\lambda_1,\ldots,\lambda_n$ of the
linear part we show that the convergence radius $\rho$ of the
conjugating transformation satisfies $\ln \rho(\lambda )\ge
-C\Gamma(\lambda)+C'$ with $\Gamma(\lambda)$ characterizing the
eigenvalues $\lambda$,  a constant $C'$ not depending on
$\lambda$ and  $C=1$.  This improves the previous results for $n\gt 1$, where the
known proofs give $C=2$.  We also recall that $C=1$ is known to be the
optimal value for $n=1$.}

\section{1}{Introduction and statement of the result}
We reconsider the classical problem of iteration of analytic
functions, previously investigated by
Schr\"oder~\dbiref{Schroeder-1871} and Siegel~\dbiref{Siegel-1942},
extending it to higher dimension.  Our aim is to improve the existing
results on the convergence radius of an analytic near the identity
transformation that conjugates the map to its linear part, thus
producing general and possibly optimal estimates.

We consider an analytic map of the complex space $\complessi^n$ into
itself that leaves the origin fixed, so that it may be written as
$$
x' = \Lambdamat x + v_1(x) + v_2(x) + \ldots\ ,\quad x\in\complessi^n
\formula{map.1}
$$
where $\Lambdamat$ is a $n\times n$ complex matrix and $v_s(x)$ is a
homogeneous polynomial of degree $s+1$. The series is assumed to be
convergent in a neighborhood of the origin of $\complessi^n$.  We
also assume that $\Lambdamat=\diag(\lambda_1,\ldots,\lambda_n)$ is a
diagonal complex matrix with non resonant eigenvalues, in a sense to
be clarified in a short.

The problem raised by Schr\"oder is the following: to find an analytic
near the identity coordinate transformation written as a convergent
power series
$$
y_j = x_j + \phi_{1,j}(x) + \phi_{2,j}(x) + \ldots
\ ,\quad j=1,\ldots,n\ ,
\formula{lser.7}
$$
with $\phi_{s,j}$ of degree $s+1$ which conjugates the map~\frmref{map.1}
to its linear part, namely in the new coordinates $y$ the map takes the
form
$$
y' = \Lambdamat y\ .
\formula{sielie.101}
$$

Let us write the eigenvalues of the complex matrix $\Lambdamat$ in
exponential form, namely as $\lambda_j=e^{\mu_j+i\omega_j}$, with real
$\mu_j$ and $\omega_j$.  We define the sequence $\{\beta_r\}_{r\ge 0}$
of positive real numbers as
$$
\beta_0=1\ ,\quad 
\beta_r=\min_{|k|=r+1} 
 \bigl|e^{\langle k,\mu+i\omega\rangle-\mu_j-i\omega_j}-1\bigr|\>,\ r\ge 1
\formula{sielie.50}
$$
The eigenvalues are said to be non resonant if $\beta_r\ne 0$ for $r\ge
1$.  This is enough in order to positively answer the question raised
by Schr\"oder {\corsivo in formal sense}, i.e., disregarding the
question of convergence.

We briefly recall some known results that motivate
the present work.  For a detailed discussion of the problem we refer to
Arnold's book~\dbiref{Arnold-1980}.    
The eigenvalues are said to belong to the {\corsivo Poincar\'e domain}
in case all $\mu_j$'s have the same sign (i.e., the eigenvalues
$\lambda$ are all inside or all outside the unit circle in the complex
plane).  The complement of the Poincar\'e domain is named {\corsivo
Siegel domain}.  In the latter case the problem raised by Schr\"oder
represented a paradigmatic  case of perturbation series involving small
divisors, since the sequence $\{\beta_r\}_{r\ge 0}$ above has zero as
inferior limit.  

The formal solution of the problem has been established by
Schr\"oder~\dbiref{Schroeder-1871} in 1871 for the case $n=1$, when the
Siegel domain is the unit circle, i.e., $\lambda=e^{i\omega}$.  The
problem of convergence in this case has been solved by
Siegel~\dbiref{Siegel-1942} in 1942, assuming a strong non resonance
condition of diophantine type.  Siegel's proof is worked out using the
classical method of majorants introduced by Cauchy, and makes use of a
delicate number-theoretical lemma, often called Siegel's lemma,
putting particular emphasis on the relevance of Diophantine
approximations.  This was actually the first proof of convergence of a
meaningful problem involving small divisors.

A condition of diophantine type was also used by Kolmogorov in his
proof of persistence of invariant tori under small perturbation of
integrable Hamiltonian systems, and has become a standard one in KAM
theory.  It should be noted that Kolmogorov introduced also an
iteration scheme exhibiting a fast convergence, that he described as
``similar to Newton's method'' and is often referred to as ``quadratic
method''.  The fast convergence scheme applies also to the case
of maps discussed here: see, e.g.,~\dbiref{Arnold-1980}, \S$\,$28.

In the last decades two problems have been raised (among many
others).  The first one is concerned with determining the optimal non
resonance condition; the second one with the size of the convergence
radius of the transformation to normal form.  In the case
$n=1$ both questions have been answered exploiting the geometric
renormalization approach introduced by Yoccoz
(see~\dbiref{Yoccoz-1995}, \dbiref{Yoccoz-2002} ~\dbiref{Buff-2006}).
The optimal set of rotation numbers for which an analytic
transformation to linear normal form exists has been identified with
the set of Bruno numbers.  These numbers obey a non resonance condition
weaker than the diophantine one introduced by Siegel.  

Bruno's condition in case $n\gt 1$ may be stated as follows.  Let us
introduce the sequence $\{\alpha_r\}_{r\ge 0}$ defined as
$$
\alpha_r = \min_{0\le s\le r} \beta_s\ ,\quad r\ge 0\ .
\formula{sielie.51}
$$
The latter sequence is monotonically non increasing, and if the
eigenvalues belong to Siegel's domain has zero as lower limit.  The
non resonance condition is
$$
-\sum_{k\gt 0}\frac{\ln\alpha_{2^k-1}}{2^{k}} = \bcyr\lt \infty\ .
$$
The series expansion of the transformation that solves
Schr\"oder's problem is proved to be convergent in a disk
$\Delta_{\rho}$ centered at the origin with radius $\rho(\lambda)$
satisfying 
$$
\ln \rho \ge -C \bcyr +C'\ ,
\formula{raggio}
$$
where $C'$ is a constant independent of $\lambda$.  For $n=1$ the
optimal value has been proved to be $C=1$.  However, the geometric
renormalization methods can not be extended to the case $n\gt 1$, for
which only the value $C=2$ has been found.  A proof that gives $C=1$
in the case $n=1$ using the majorant method has been given
in~\dbiref{Carletti-Marmi-2000}, but has not been extended to the case
$n\gt 1$.

In this paper we prove that the same bound with $C=1$ holds true in
any finite dimension.  We obtain this result by implementing a
representation of the map with the method of Lie series and Lie
transforms, and producing convergence estimates in the spirit of
Cauchy's majorants method.  The formal implementation of the
representation of maps may be found in~\dbiref{Giorgilli-2013}.

This paper extends to the case of maps a previous similar result for
the case of the Poincar\'e-Siegel theoretical center problem, where
the problem of linearization of an analytic system of differential
equations in the neighborhood of a fixed point is considered.  We
stress that the interest of the method is not limited to the cases
mentioned here.  E.g.: for applications to KAM theory
see~\dbiref{Giorgilli-1997.4}, \dbiref{Giorgilli-1999}
and~\dbiref{Gio-Loc-San-2013};  a proof of Lyapounov's theorem on the
existence of periodic orbits in the neighborhood of the equilibrium
is given in~\dbiref{Giorgilli-2012};  extensions of Lyapounov's
theorem may be found in~\dbiref{Giorgilli-2001}
and~\dbiref{Gio-Mur-2006}.

We come now to a formal statement of our result.  We assume the
following

\noindent
{\bf Condition $\tauvet$:\enspace}The sequence $\alpha_r$ above satisfies
$$
-\sum_{r\ge 1} \frac{\ln \alpha_r}{r(r+1)} = \Gamma \lt \infty\ .
\formula{Ctau}
$$

\theorem{mainth}
{Consider the map~\frmref{map.1} and assume that the eigenvalues of
$\Lambdamat$ are non-resonant and satisfy condition $\tauvet$. Then there
exists a near to identity coordinate transformation $y=x+\psi(x)$,
with $\psi$ analytic at least in the polydisk of radius $B^{-1}e^{-\Gamma}$,
where $B>0$ is a universal constant, which transforms the map into the
normal form $y'=\Lambdamat y$.}
\endclaim

Our condition~$\tauvet$ is equivalent to Bruno's one.  We have indeed
$$
\Gamma\le\bcyr\le 2\Gamma\ .
$$
However, our formulation of condition $\tauvet$ comes out naturally
from our analysis of accumulation of small divisors and turns out to
be the key that allows us to find the estimate with $C=1$.  For a
brief discussion of the relations between the two conditions
see~\dbiref{Gio-Mar-2010}.

The paper is organized as follows.  In section~\secref{2} we include
the essential information on the formal algorithm, referring
to~\dbiref{Giorgilli-2013} for details.  In section~\secref{3} we give
the technical estimates that lead to the proof of convergence of the
formal algorithm.  In particular we include a complete discussion of
the mechanism of accumulation of small divisors which allows us to
have an accurate control.  We remark in particular that we do not need
to use the Siegel lemma, because the latter controls combinations
of small divisors that do not occur in our scheme.  The estimates of
section~\secref{3} are used in section~\secref{4} in order to complete
the proof of the main theorem.  A technical appendix follows.

\section{2}{Formal algorithm}
We represent a map of type~\frmref{map.1} using Lie transforms.  A
detailed exposition of the representation algorithm and of the method
of reduction to normal form is given in~\dbiref{Giorgilli-2013}, where
the problem is treated at a formal level.  We refer to that paper for
details concerning the formal setting, while in this paper we pay
particular attention to the problem of convergence.  Here we just
include the definitions that we shall need later, so that also the
notations will be fixed.  Some relevant lemmas are reported in
appendix~\appref{A}.

We start with the definition of Lie series and Lie transform.  Let
$X_s(x)$ be a vector field on $\complessi^n$ whose components are
homogeneous polynomials of degree $s+1$.  We shall say that $X_s(x)$
is of {\corsivo order} $s$, as indicated by the label.  Moreover, in
the following we shall denote by $X_{s,j}$ the $j$-th component of
the vector field $X_s\,$.  The {\corsivo Lie series} operator is
defined as
$$
\exp(\lie{X_s}) = \sum_{j\ge 0} \frac{1}{j!} \lie{X_s}^j
\formula{sielie.100}
$$
where $\lie{X_s}$ is the Lie derivative with respect to the vector
field $X_s$.

Let now $X=\{X_j\}_{j\ge 1}$ be a sequence of polynomial vector fields
of degree $j+1$.  The {\corsivo Lie transform} operator $T_X$ is
defined as
$$
T_{X} = \sum_{s\ge 0} E^{X}_s\ ,
\formula{trslie.1}
$$
where the sequence $E^{X}_s$ of linear operators is recursively
defined as
$$
E^{X}_0 = \uno\ ,\quad 
E^{X}_s = \sum_{j=1}^{s} \frac{j}{s}\lie{X_j} E^{X}_{s-j}\ .
\formula{trslie.2}
$$
The superscript in $E^{X}$ is introduced in order to specify which
sequence of vector fields is intended.  By letting the sequence to
have only one vector field different from zero, e.g.,
$X=\{0,\ldots,0,X_k,0,\ldots\}$ it is easily seen that one gets
$T_X=\exp\bigl(\lie{X_k}\bigr)$.

\subsection{2.2}{Representation and conjugation of maps}
We recall the representation of maps introduced
in~\dbiref{Giorgilli-2013} together with some formal results that we
are going to use here.  Let $\Lambdamat=e^{\Amat}$, i.e., in our case,
$\Amat=\diag(\mu_1+i\omega_1,\ldots,\mu_n+i\omega_n)$ with
$\lambda_j=e^{\mu_j+i\omega_j}$.  Remark that we may express the
linear part of the map as a Lie series by introducing the exponential
operator $\Rmat=\exp\bigl(\lie{\Amat x}\bigr)$.  The action of the
operator $\Rmat$ on a function $f$ or on a vector field $V$ is easily
calculated as
$$
\bigl(\Rmat f\bigr)(x) = f(\Lambdamat x)\ ,\quad
 \bigl(\Rmat V\bigr)(x) = \Lambdamat^{-1} V(\Lambdamat x)\ .
\formula{sielie.112}
$$

The first result is concerned with the representation of the
map~\frmref{map.1} using a Lie transform.

\lemma{map.13}{There exist generating sequences of
vector fields $V= \bigl\{V_s(x)\bigr\}_{s\geq1}$ and $W
=\bigl\{W_s(x)\bigr\}_{s\geq1}$ with $W_s =\Rmat V_s$ such that
the map~\frmref{map.1} is represented as
$$
x' = \Rmat\circ T_{V} x\quad {\rm and}\quad
x' = T_{W}\circ\Rmat\, x
\formula{map.5}
$$
}\endclaim

The second result is concerned with the composition of Lie
transforms. 

\lemma{trslie.11}{Let $X,\,Y$ be generating sequences.  Then one
has $T_X\circ T_Y = T_Z$ where $Z$ is the generating sequence
recursively defined as
$$
Z_1 = X_1 + Y_1\ ,\quad
Z_s = X_s + Y_s + \sum_{j=1}^{s-1} \frac{j}{s} E^{X}_{s-j} Y_j\ .
\formula{trslie.12}
$$
}\endclaim

\noindent
The latter formula reminds the well known Baker-Campbell-Hausdorff
composition of exponentials.  The difference is that the result is
expressed as a Lie transform instead of an exponential, which makes
the formula more effective for our purposes.

The third result gives the algorithm that we shall use in order to
conjugate the map~\frmref{map.1} to its linear part.  We formulate it
in a more general form, looking for conjugation between maps.  Let two
maps
$$
x'=T_{W}\circ\Rmat\, x\ ,\quad y'=T_{Z}\circ \Rmat\, y
\formula{map.6}
$$ 
be given, where 
$W=\bigl\{W_s\bigr\}_{s\geq1}$, $Z=\bigl\{Z_s\bigr\}_{s\geq1}$  are generating
sequences.  We say that the maps are conjugated up to order $r$ in
case there
exists a finite generating sequence $X=\bigl\{X_1,\ldots,X_r\bigr\}$ such that
the transformation $y=S^{(r)}_{X}x$ makes the difference between the maps
to be of order higher than $r$, i.e.,
$$
S^{(r)}_{X} x'\Big|_{x' = T_{W}\circ\Rmat\, x} 
 - T_{Z}\circ\Rmat\, y\Big|_{y = S^{(r)}_{X} x} = \Oscr(r+1)\ ,
$$
where $S^{(r)}_{X}
=\exp\bigl(\lie{X_r}\bigr)\circ\ldots\circ\exp\bigl(\lie{X_1}\bigr)$.
Suppose that we have $W_1=Z_1,\ldots,W_{r}=Z_{r}$.  Then the maps
are  conjugated up to order $r$, since one has
$T_{W} x - T_{Z} x = \Oscr(r+1)$.

\lemma{map.20}{Let the generating sequences of the
maps~\frmref{map.6} coincide up to order $r-1$ and let $X_r$ be a
vector field of order $r$ generating the near the identity
transformation $y=\exp\bigl(\lie{X_r} x\bigr)$.  Then the maps are
conjugated up to order $r$ if
$$
T_{Z} =  
   \exp\bigl(\lie{X_r}\bigr) \circ T_{W}
    \circ\exp\bigl(\lie{-\Rmat X_r}\bigr)\ .
\formula{map.21a}
$$
The vector field $X_r$ must satisfy the equation
$$
\Dmat X_r  = W_r - Z_r\ ,\quad \Dmat = \Rmat-\uno\ .
\formula{map.21b}
$$}\endclaim

\subsection{2.3}{Construction of the normal form}
Following Schr\"oder and Siegel, we want to conjugate the
map~\frmref{map.1} to its linear part.  That is: writing the map as 
$$
x' =T_{W^{(0)}}\circ\Rmat x
\formula{sielie.102}
$$
with a known sequence $W^{(0)}$ of vector fields, we want to reduce it
to the linear normal form
$$
x' = \Rmat x\ .
\formula{sielie.103}
$$
To this end we look for a generating sequence $\{X_r\}_{r\ge 1}$ of
vector fields and a corresponding sequence $\{W^{(r)}\}_{r\ge 1}$
satisfying $W^{(r)}_1=\ldots=W^{(r)}_r=0$.  We emphasize that the map
$x' =T_{W^{(r)}}\circ\Rmat x$ is conjugated to~\frmref{sielie.103} up
to order $r$.  We say that $W^{(r)}$ is in normal form up to order
$r$.

According to~\frmref{map.21b} we should solve for $X_r$ the equation 
$$
\Dmat X_r = W^{(r-1)}_r\ ,\quad \Dmat = \Rmat - \uno\ .
\formula{sielie.107}
$$
The operator $\Dmat$ is diagonal on the basis of
monomials $x^k\evet_j = x_1^{k_1}\cdot\ldots\cdot x_n^{k_n}\evet_j$,
where $(\evet_1,\ldots,\evet_n)$ is the canonical basis of
$\complessi^n$.  For we have
$$
\Dmat\, x^k\evet_j 
 = \bigl(e^{\langle k,\mu+i\omega\rangle - \mu_j-i\omega_j}-1\bigr)\,
    x^k\evet_j \ .
\formula{sielie.108}
$$
Thus, provided the eigenvalues of $\Dmat$ are not zero, the vector
field $X_r$ is determined as
$$
X_r = \sum_{j=1}^{n} \evet_j \sum_{k} 
 \frac{w_{j,k}}{e^{\langle k,\mu+i\omega\rangle - \mu_j-i\omega_j}-1}
 x^k\ ,
$$
where $w_{j,k}$ are the coefficients of $W^{(r-1)}_r$.

Next, we need an explicit form for the transformed sequence $W^{(r)}$.
We use the conjugation formula~\frmref{map.21a} replacing $W$ and $Z$
with $W^{(r-1)}$ and $W^{(r)}$, respectively.  It is convenient to
introduce an auxiliary vector field $V^{(r)}$ and to split the formula
as
$$
\eqalign{
T_{V^{(r)}} &=
   T_{W^{(r-1)}}
    \circ \exp\bigl(\lie{-\Rmat X_r}\bigr)\ .
\cr
T_{W^{(r)}} &= 
   \exp\bigl(\lie{X_r}\bigr) \circ T_{V^{(r)}}\ ,
\cr
}
\formula{sielie.5}
$$
A more explicit form comes from lemma~\lemref{trslie.11}, recalling
that $T_{X}=\exp(\lie{X_r})$ if one considers the generating sequence
$X=\{0,..,0,X_r,0,\ldots\}$.  The auxiliary vector field $V^{(r)}$ is
determined as
$$
\vcenter{\openup1\jot\halign{
\hfil$\displaystyle{#}$
&$\displaystyle{#}$\hfil
&\quad\quad{\rm for\ }$\displaystyle{#}$\hfil
\cr
V^{(r)}_r 
&=W^{(r-1)}_r - \Rmat X_r\ ,
\cr
V^{(r)}_s
&= W^{(r-1)}_s - \frac{r}{s} E^{(r-1)}_{s-r} \Rmat X_r & s\gt r\>,
\cr
}}
\formula{sielie.8}
$$
where we use the simplified notation $E^{(r-1)}$ in place of
$E^{(W^{(r-1)})}$.  
Having so determined the sequence $V^{(r)}$ we calculate the
transformed sequence $W^{(r)}$ as
$$
\vcenter{\openup1\jot\halign{
\hfil$\displaystyle{#}$
&$\displaystyle{#}$\hfil
&\quad{\rm for\ }$\displaystyle{#}$\hfil
\cr
W^{(r)}_r
&= V^{(r)}_r + X_r \ ;
\cr
W^{(r)}_{s}
&= V^{(r)}_{s} + \frac{1}{s} \sum_{k=1}^{\lfloor s/r\rfloor-1} 
 { \frac{s-kr}{k!}} \lie{X_r}^{k} V^{(r)}_{s-kr}
\ ,\quad s\gt r\>.
\cr
}}
\formula{sielie.6}
$$
Here a remark is in order. The formul{\ae} above define the sequences
$V^{(r)}_r,V^{(r)}_{r+1},\ldots$ and
$W^{(r)}_{r+1},W^{(r)}_{r+2},\ldots$ starting with terms of order $r$
and $r+1$, respectively.  This is fully natural, because all terms of
lower order vanish due to $W^{(r)}$ being in normal form up to order
$r$.  Moreover we emphasize that in view of~\frmref{map.21b} we should
determine $X_r$ by solving the equation $\Dmat X_r=W^{(r-1)}_r$, since
we want $W^{(r)}_r=0$.

\section{3}{Quantitative estimates}
Our aim now is to complete the formal algorithm of the previous
section with quantitative estimates that will lead to the proof of
convergence of the transformation to normal form.  The main result of
this section is the iteration lemma of section~\sbsref{3.3} below.
However, we must anticipate a few technical tools.

\subsection{3.1}{Norms on vector fields and generalized Cauchy estimates}
For a homogeneous polynomial $f(x)=\sum_{|k|=s}f_kx^k$ (using
multiindex notation and with $|k|=|k_1|+\ldots+|k_n|$) with complex
coefficients $f_k$ and for a homogeneous polynomial vector field
$X_s=(X_{s,1},\ldots,X_{s,n})$ we use the {\it polynomial norm}
$$
\bignorma{f} = \sum_{k} |f_k|\ ,\quad
 \bignorma{X_s} = \sum_{j=1}^{n}\, \bignorma{X_{s,j}}\ .
\formula{sielie.104}
$$

The following lemma allows us to control the norms of Lie derivatives
of functions and vector fields.
 
\lemma{sielie.105}{Let $X_r$ be a homogeneous polynomial vector field
of degree $r+1$.  Let $f_s$ and $v_s$ be a homogeneous polynomial and
vector field, respectively, of degree $s+1$.  Then we have
$$
\bigl\|\lie{X_r} f_s\bigr\| \le 
 (s+1)\, \|X_r\|\, \|f_s\|\quad {\rm and}\quad
  \bignorma{\liec{X_r} v_s}
 \le (r+s+2)\,
  \|X_r\|\, \|v_s\|\ .
\formula{sielie.106}
$$}\endclaim

\proof
Write $f_s=\sum_{|k|=s+1} b_k x^k$ with complex coefficients
$b_k$. Similarly, write the $j$-th component of the vector field
$X_r$ as $X_{r,j}=\sum_{|k'|=r+1} c_{j,k'} x^{k'}$.  Recalling that
$\lie{X_r}f_s=\sum_{j=1}^{n}X_{r,j}\derpar{f_s}{x_j}$  we have
$$
\lie{X_r}f_s = \sum_{j=1}^{n} \sum_{k,k'}
 \frac{c_{j,k'}k_jb_k}{x_j} x^{k+k'}\ .
$$
Thus in view of $|k_j|\leq s+1$ we have
$$
\bignorma{\lie{X_r}f_s} 
 \le 
   (s+1)\sum_{j=1}^{n}\sum_{k'} |c_{j,k'}| \sum_{k} |b_k| 
   = (s+1)\ordnorma{X_r}\, \ordnorma{f_s}\ ,
$$
namely the first of~\frmref{sielie.106}.  In order to prove the second
inequality recall that the $j$-th component of the Lie derivative of
the vector field $v_s$ is
$$
\bigl(\lie{X_r} v_s\bigr)_j =
 \sum_{l=1}^{n}\left(X_{r,l}\derpar{v_{s,j}}{x_l} 
  - v_{s,l}\derpar{X_{r,j}}{x_l}\right)\ .
$$
Then using the first of~\frmref{sielie.106} we have
$$
\eqalign{
\biggnorma{\sum_{l=1}^{n}\left(X_{r,l}\derpar{v_{s,j}}{x_l} 
  - v_l\derpar{X_{r,j}}{x_l}\right)}
\le  (s+1) \ordnorma{X_r}\, \ordnorma{v_{s,j}}
    +(r+1) \ordnorma{v_s}\, \ordnorma{X_{r,j}}\ , 
\cr
}
$$
which readily gives the wanted inequality in view of the
definition~\frmref{sielie.104} of the polynomial norm of a vector
field.\endproof

\lemma{sielie.109}{Let $V_r$ be a homogeneous polynomial vector field
of degree $r+1$.  Then the solution $X_r$ of the equation $\Dmat X_r
= V_r$ satisfies
$$
\bignorma{X_r} \le \frac{1}{\alpha_r}\bignorma{V_r}\ ,\quad
 \bignorma{\Rmat X_r} \le \frac{1+\alpha_r}{\alpha_r}\bignorma{V_r}\ ,
\formula{sielie.110}
$$
with the sequence $\alpha_r$ defined by~\frmref{sielie.51}
}\endclaim

\proof
The first inequality is a straightforward consequence of the
definition~\frmref{sielie.104} of the norm and of the sequence
$\alpha_r$ in terms of $\beta_r$ defined by~\frmref{sielie.50}.  For
if $v_{j,k}$ are the coefficients of $V_r$ then the coefficients of
$X_r$ are bounded by $|v_{j,k}|/\beta_r\le |v_{j,k}|/\alpha_r$.  The
second inequality follows from $\Rmat X_r = X_r+V_r$, which gives the
stated inequality.
\endproof

\subsection{3.2}{Accumulation of small divisors}
Lemma~\lemref{sielie.109} shows that solving the equation for every
vector field of the generating sequence introduces a divisor
$\alpha_r$.  Such divisors do accumulate, and our aim now is to
analyze in detail the process of accumulation.  It will be
convenient to introduce a further sequence $\{\sigma_r\}_{r\ge 0}$
defined as
$$
\sigma_0 = 1\ ,\quad\sigma_r = \frac{\alpha_r}{r^2}\>,\ r\ge 1\>,
\formula{sielie.111}
$$
which will play a major role in the rest of the proof.  The extra
factor $1/r^2$ can be interpreted as due to the generalized Cauchy
estimates for derivatives, that are also a source of divergence in
perturbation processes.  The quantities $\sigma_r$ are the actual
small divisors that we must deal with.  Here we follow the scheme
presented in~\dbiref{Gio-Mar-2010}, with a few {\corsivo variazioni}.
However, since this is a crucial part of the proof we include it in a
detailed and self contained form.

The guiding remark is that the small divisors $\sigma_r$ propagate and
accumulate through the formal construction due to the use of the
recursive formul{\ae}~\frmref{sielie.8} and~\frmref{sielie.6}.  E.g.,
the expression $\lie{X_r}V_s^{(r)}$ will contain the product of the
divisors contained in both $X_r$ and $V_s^{(r)}$, with no extra
divisors generated by the Lie derivative.  The explicit constructive
form of the algorithm allows us to have a quite precise control on the
accumulation process.  It is an easy remark that unfolding the
recursive formul{\ae}~\frmref{sielie.8} and~\frmref{sielie.6} will
produce an estimate of $\ordnorma{W^{(r)}_{s}}$ as a sum of many terms
every one of which contains as denominator a product of $q$ divisors
of the form $\sigma_{j_1}\cdots\sigma_{j_q}$, with some indexes
$j_1,\ldots,j_q$ and some $q$ to be found.  This is what we call the
{\corsivo accumulation of small divisors}, and the problem is to
identify the worst product among them.  The key of our argument is to
focus our attention on the indexes rather than on the actual values of
the divisors.

We call $I=\{j_1,\ldots,j_s\}$ with non negative integers
$j_1,\ldots,j_s$ a {\sl set of indexes}.  We introduce a partial
ordering as follows.  Let $I=\{j_1,\ldots,j_s\}$ and
$I'=\{j'_1,\ldots,j'_s\}$ be two sets of indexes with the same number
$s$ of elements.  We say that $I\ilt I'$ in case there is a
permutation of the indexes such that the relation $j_m\le j'_m$ holds
true for $m=1,\ldots,s\,$.  If two sets of indexes contain a different
number of elements we pad the shorter one with zeros and use the same
definition.
We also define the special sets of indexes
$$
I^*_s = \Bigl(\Bigl\lfloor\frac{s}{s}\Bigr\rfloor, 
               \Bigl\lfloor\frac{s}{s-1}\Bigr\rfloor,\ldots, 
                \Bigl\lfloor\frac{s}{2}\Bigr\rfloor \Bigr)\ .
\formula{sielie.53}
$$

\lemma{nrmlie.44}{For the sets of indexes $I_s^*=\{j_1,\ldots,j_s\}$ the
following statements hold true:
\item{(i)}the maximal index is  $j_{\rm
max}=\bigl\lfloor\frac{s}{2}\bigr\rfloor\,$;
\item{(ii)}for every $k\in\{1,\ldots,j_{\max}\}$ the index $k$
appears exactly
$\bigl\lfloor\frac{s}{k}\bigr\rfloor -\bigl\lfloor\frac{s}{k+1}\bigr\rfloor$
times\/;
\item{(iii)}for $0\lt r\le s$ one has
$$
\bigl(\{r\}\cup I^*_r\cup I^*_s\bigr) \ilt I^*_{r+s}\ .
$$
}\endclaim

\proof
The claim~(i) is a trivial consequence of the definition.

\noindent
(ii) For each fixed value of $s>0$ and $1\le k\le\lfloor
s/2\rfloor\,$, we have to determine the cardinality of the set
$M_{k,s}=\{m\in\naturali:\ 2\le m\le s\,,\ \lfloor s/m\rfloor=k\}$.
For this purpose, we use the trivial inequalities
$$
\biggl\lfloor\frac{s}{\lfloor s/k\rfloor}\biggr\rfloor \geq k
\quad\hbox{and}\quad
\biggl\lfloor\frac{s}{\lfloor s/k\rfloor+1}\biggr\rfloor \lt k\ .
$$
After having rewritten the same relations with $k+1$ in place of
$k\,$, one immediately realizes that a index $m\in M_{k,s}$ {\sl if
and only if} $m\le\lfloor s/k\rfloor$ {\sl and} $m\ge \lfloor
s/(k+1)\rfloor+1\,$, therefore
$\#M_{k,s}=\bigl\lfloor\frac{s}{k}\bigr\rfloor
-\bigl\lfloor\frac{s}{k+1}\bigr\rfloor\,$.

\noindent
(iii)~Since $r\le s\,$, the definition in~\frmref{sielie.53} implies
that neither $\{r\}\cup I^*_r\cup I^*_s$ nor $I^*_{r+s}\,$ can include
any index exceeding $\bigl\lfloor(r+s)/2\bigr\rfloor\,$.  Thus, let us
define some finite sequences of non-negative integers as follows:
$$
\vcenter{\openup1\jot\halign{
 \hbox {\hfil $\displaystyle {#}$}
&\hbox {\hfil $\displaystyle {#}$\hfil}
&\hbox {$\displaystyle {#}$\hfil}
&\hbox to 2 ex{\hfil$\displaystyle {#}$\hfil}
&\hbox {\hfil $\displaystyle {#}$}
&\hbox {\hfil $\displaystyle {#}$\hfil}
&\hbox {$\displaystyle {#}$\hfil}\cr
R_k &= &\#\bigl\{j\in I^*_{r}\>:\> j\le k\bigr\}\ ,
& &S_k &= &\#\bigl\{j\in I^*_{s}\>:\> j\le k\bigr\}\ ,
\cr
M_k &= &\#\bigl\{j\in \{r\}\cup I^*_r\cup I^*_s\>:\> j\le k\bigr\}\ ,
& &N_k &= &\#\bigl\{j\in I^*_{r+s}\>:\> j\le k\bigr\}\ ,
\cr
}}
$$
where $1\leq k \leq\lfloor(r+s)/2\rfloor\,$. When $k\lt r\,$, the
property~(ii) of the present lemma allows us to
write
$$
R_k = r - \Bigl\lfloor \frac{r}{k+1}\Bigr\rfloor\ ,\quad
S_k = s - \Bigl\lfloor \frac{s}{k+1}\Bigr\rfloor\ ,\quad
N_k = r+s - \Bigl\lfloor \frac{r+s}{k+1}\Bigr\rfloor\ ;
$$
using the elementary estimate $\lfloor x\rfloor + \lfloor
y\rfloor \leq \lfloor x+y\rfloor\,$, from the equations above it
follows that $M_k \geq N_k$ for $1\le k<r\,$.  In the remaining
cases, i.e., when $r\le k\le \lfloor(r+s)/2\rfloor\,$, we have that
$$
R_k = r - 1\ ,\quad
S_k = s - \Bigl\lfloor \frac{s}{k+1}\Bigr\rfloor\ ,\quad
N_k = r+s - \Bigl\lfloor \frac{r+s}{k+1}\Bigr\rfloor\ ;
$$
therefore, $M_k=1+R_k+S_k\ge N_k\,$.
\noindent
Since we have just shown that $M_k \geq N_k$ $\forall\ 1\le
k\le \lfloor(r+s)/2\rfloor\,$, it is now an easy matter to complete
the proof. Let us first imagine to have reordered both the set of
indexes $\{r\}\cup I^*_r\cup I^*_s$ and $I^*_{r+s}$ in increasing order;
moreover, let us recall that $\#\big(\{r\}\cup I^*_r\cup I^*_s\big)=\#
I^*_{r+s}=r+s-1\,$, because of the definition in~\frmref{sielie.53}.
Thus, since $M_1\ge N_1\,$, every element equal to $1$ in $\{r\}\cup
I^*_r\cup I^*_s$ has a corresponding index in $I^*_{r+s}$ the value of
which is at least $1\,$. Analogously, since $M_2\ge N_2\,$, every
index $2$ in $\{r\}\cup I^*_r\cup I^*_s$ has a corresponding index in
$I^*_{r+s}$ which is at least $2\,$, and so on up to $k
=\lfloor(r+s)/2\rfloor\,$. We conclude that $\{r\}\cup
I^*_r\cup I^*_s\ilt I^*_{r+s}\,$.\endproof

We come now to identify the sets of indexes that describe the allowed
combinations of small divisors.  These are the sets
$$
\Jscr_{r,s} = \bigl\{I=\{j_1,\ldots,j_{s-1}\}\>:\>
 j_m\in\{0,\ldots,\min(r,s/2)\}\,,\> I\ilt I^*_s\bigr\}\ .
\formula{nrmlie.20}
$$
We shall refer to condition $I\ilt I^*_s$ as the {\corsivo
selection rule $\Smat$}.  The relevant properties are stated by the
following

\lemma{nrmlie.23}{For the sets of indexes $\Jscr_{r,s}$ the following
statements hold true:
\item{(i)}$\Jscr_{r-1,s}\subset\Jscr_{r,s}\,$;
\item{(ii)}if $I\in\Jscr_{r-1,r}$ and $I'\in\Jscr_{r,s}$ then
we have $\bigr(\{r\} \cup I \cup I'\bigl)\in\Jscr_{r,r+s}\,$.}\endclaim

\remark
Property~(ii) plays a major role in controlling the accumulation
of small divisors, since it gives us a control of how the indexes are
accumulated.  For it reflects the fact that a Lie derivative, e.g.,
$\lie{X_r}V_s^{(r)}$ contains the union of the divisors in $V_s^{(r)}$ and in
$X_r$, taking also into account that $X_r$ contains an
extra divisor in view of lemma~\lemref{sielie.109}.

\prooftx{of lemma~\lemref{nrmlie.23}}
(i)~is immediately checked in view of the definition of $\Jscr_{r,s}$.

\noindent
(ii)~We have $\#\bigl(\{r\}\cup\, I\cup\, I'\bigr) = 1+\#(I)+\#(I') =
1+r-1+s-1 = r+s-1\,$.  For $j\in\{r\}\cup\, I\cup\, I'\,$ we also have
$1\le j\le r$ because this is true for all $j\in I$ and for all $j\in
I'$, and we just add an extra index $r$.  Since $r\lt s$, we also have
$r =\min\bigl(r,(r+s)/2\bigr)$, as required.  Coming to the selection
rule $\Smat$ we remark that $\{r\}\cup I\cup I'\ilt\{r\}\cup I^*_r\cup
I^*_s$ readily follows from $I\ilt I^*_r$ and $I'\ilt I^*_s$, which
are true in view of $I\in\Jscr_{r-1,r}$ and $I'\in\Jscr_{r,s}$, so
that the claim follows from property~(iii) of
lemma~\lemref{nrmlie.44}.\endproof

We come now to consider the accumulation of small divisors.  Recall
the definition~\frmref{sielie.111} of the sequence $\sigma_r$.
We associate to the sets of indexes $\Jscr_{r,s}$ the sequence of
positive real numbers $T_{r,s}$ defined as
$$
T_{0,s} = 1\ ,\quad
T_{r,s} = \max_{I\in\Jscr_{r,s}} \prod_{j\in I} \frac{1}{\sigma_j}
 \ ,\quad  0\lt r\leq s\ .
\formula{nrmlie.26}
$$

\lemma{nrmlie.27}{The sequence $T_{r,s}$ satisfies the following
properties for $1\leq r\leq s\,$:
\item{(i)}$T_{r-1,s} \le T_{r,s}\,$;
\item{(ii)}$\frac{1}{\sigma_r} T_{r-1,r} T_{r,s} \le
T_{r,r+s}\,$.
}\endclaim

\proof
(i)~From property~(i) of lemma~\lemref{nrmlie.23} we readily get
$$
T_{r-1,s}
= \max_{I\in\Jscr_{r-1,s}} \prod_{j\in I}\frac{1}{\sigma_j} \le
\max_{I\in\Jscr_{r,s}} \prod_{j\in I}\frac{1}{\sigma_j}\ ,
$$
since the maximum is evaluated over a larger set of indexes.

\noindent
(ii)~Compute
$$
\eqalign{
\frac{1}{\sigma_r} T_{r-1,r} T_{r,s}
&= \frac{1}{\sigma_r}
    \max_{I\in\Jscr_{r-1,r}} \prod_{j\in I}\frac{1}{\sigma_j}
     \max_{I'\in\Jscr_{r,s}} \prod_{j'\in I'}\frac{1}{\sigma_{j'}}
\cr
& = \max_{I\in\Jscr_{r-1,r}} \max_{I'\in\Jscr_{r,s}}
     \prod_{j\in\{r\}\cup\,I\cup\,I'}\frac{1}{\sigma_{j}}
\cr
& \le \max_{I\in\Jscr_{r,r+s}} \prod_{j\in I}\frac{1}{\sigma_j}
= T_{r,r+s}\ ,
\cr
}
$$
where in the inequality of the last line property~(ii) of
lemma~\lemref{nrmlie.23} has been used.\endproof

The final estimate uses condition~$\tauvet$ and the definition~\frmref{sielie.111} of the
sequence $\sigma_r\,$.

\lemma{nrmlie.42}{Let $\lambda$ satisfy condition~$\tauvet$.  Then
the sequence $T_{r,s}$ is bounded by
$$
T_{r,s} \le \gamma^s e^{s\Gamma}\ ,\quad
\frac{1}{\sigma_s} T_{r,s} \le \gamma^s e^{s\Gamma}
$$
with some positive constant $\gamma$ not depending on $\lambda$.}\endclaim

\proof
In view of $\sigma_s\le 1$ it is enough to prove the second
inequality.  We use the property~(ii) of lemma~\lemref{nrmlie.23} and
the selection rule~$\Smat$.  We readily get
$$
\frac{1}{\sigma_s} T_{r,s}
= \frac{1}{\sigma_s}
   \max_{I\in\Jscr_{r,s}} \prod_{j\in I} \frac{1}{\sigma_j}
\le \max_{I\in\Jscr_{r,s}} \prod_{j\in\{s\}\cup I} \frac{1}{\sigma_j}
\le \prod_{j\in\{s\}\cup I^*_s} \frac{1}{\sigma_j}\ .
$$
By property~(ii) of lemma~\lemref{nrmlie.44} the latter product is
evaluated as
$$
\prod_{j\in\{s\}\cup I^*_s} \frac{1}{\sigma_j} =
 \left[
  \sigma_1^{q_1}\cdot\ldots\cdot
   \sigma_{\scriptscriptstyle \lfloor s/2\rfloor}^{q_{\lfloor s/2\rfloor}}
    \sigma_s^{q_s}
 \right]^{-1}\ ,
$$
where $q_k
=\bigl\lfloor\frac{s}{k}\bigr\rfloor
-\bigl\lfloor\frac{s}{k+1}\bigr\rfloor$ is the number of indexes in
$I_s^*$ which are equal to $k$.  In view of the
definition~\frmref{sielie.111} of the sequence $\sigma_s$ we have
$$
\ln \frac{1}{\sigma_r} T_{r,s}
\le  -\sum_{k=1}^{s}
   \Bigl(\Bigl\lfloor\frac{s}{k}\Bigr\rfloor
-\Bigl\lfloor\frac{s}{k+1}\Bigr\rfloor\Bigr) (\ln\alpha_k - 2\ln k)
\le - s \sum_{k\ge 1} \frac{\ln\alpha_k - 2\ln k}{k(k+1)}
= s(\Gamma + a)
\ ,
$$
where $a=\sum_{k\ge 1} \frac{2\ln k}{k(k+1)}\lt\infty$ is clearly
independent of $\lambda$.  The claim
follows by just setting $\gamma=e^a\,$.\endproof

\subsection{3.3}{Iteration lemma}
This is the main lemma which allows us to control the norms of the
sequence of vector fields $\{X_r\}_{r\gt 1}\,$.

\lemma{sielie.68}{Assume that the sequence $W^{(0)}$ of vector fields
satisfies $\ordnorma{W^{(0)}_s}\le \frac{C_0^{s-1}A}{s}$ with some
constants $A\gt 0$ and $C_0\ge 0$.  Then the sequence of vector fields
$\{X_r\}_{r\ge 1}$ that for every $r$ give the normal form $W^{(r)}$
satisfies the following estimates: there exists a bounded
monotonically non decreasing sequence $\{C_r\}_{r\geq 1}$ of positive
constants, with $C_r\to C_{\infty}\lt\infty$ for $r\to\infty$, such
that we have
$$
\bignorma{X_r} \le T_{r-1,r} \frac{C_{r-1}^{r-1}A}{r\alpha_r}
 \ ,\quad
  \bignorma{W^{(r)}_s} \le T_{r,s} \frac{C_r^{s-1} A}{s}\ .
\formula{sielie.71}
$$
The sequence may be recursively defined as
$$
C_1 = 2C_0+16A\ ,\quad
 C_r = \Bigl(1+\frac{1}{r^2}\Bigr)^{1/r}\Bigl(1+\frac{1}{r}\Bigr)^{1/r} C_{r-1}\ ,
\formula{sielie.72}
$$
so that one has $C_r\gt 16A$.
}\endclaim

\proof
The proof proceeds by induction.  For $r=0$ the inequality for
$W^{(0)}_s$ is nothing but the initial hypothesis, recalling that by
definition we have $T_{0,s}=1$.  From this, the corresponding
inequality for $X_1$ immediately follows by lemma~\lemref{sielie.109}.
The induction requires two main steps, namely: (a)~estimating
$\ordnorma{V^{(r)}_s}$ as defined by~\frmref{sielie.8}, and
(b)~estimating $\ordnorma{W^{(r)}_s}$ as defined by~\frmref{sielie.6}.

\noindent
The estimate for $X_r$ is determined by solving~\frmref{sielie.107}.
In view of the induction hypothesis~\frmref{sielie.71} we rewrite the
claim of lemma~\lemref{sielie.109} as
$$
\bignorma{X_r} \le 
 \frac{r^2 T_{r-1,r}}{\alpha_r} \cdot
  \frac{C_{r-1}^{r-1} A}{r^3}\ ,\quad
\bignorma{\Rmat X_r} \le 
 \frac{T_{r-1,r}}{\alpha_r} \cdot
  \frac{(1+\alpha_r)C_{r-1}^{r-1} A}{r}\ .\quad
\formula{sielie.26}
$$

\noindent
In order to estimate $\ordnorma{V^{(r)}_s}$ we first prove that for
$s\gt r$ we have
\formdef{sielie.66}
\formdef{sielie.62}
$$
\eqalignno{
\bignorma{E_{s-1}^{(1)} \Rmat X_1} 
& \le T_{1,s} C_0^{s-1} A \,\frac{8A}{C_0} 
   \left(1+\frac{4A}{C_0}\right)^{s-2}\ ,
&\frmref{sielie.66}
\cr
\bignorma{E_{s-r}^{(r-1)} \Rmat X_r}
&\lt T_{r,s} \frac{C_{r-1}^{s-1} A}{r^2}
  \left(1+\frac{1}{r}\right)^{\frac{s}{r}-2}\quad {\rm for}\ r\gt 1\>.
&\frmref{sielie.62}
\cr
}
$$
The proof of the latter estimates is based on the general inequality
$$
\bignorma{E_{s-r}^{(r-1)} \Rmat X_r}
 \le \!\!\!\sum_{k=1}^{\lfloor s/r\rfloor-1}\!\!\! \Theta(r,s\!-\!r,k)\,
  T_{r-1,j_1}\cdots T_{r-1,j_k}\, \frac{T_{r-1,r}}{\alpha_r}
  \cdot \frac{(1\!+\!\alpha_r)C_{r-1}^{s-k-1}A^{k+1}}{r}
\formula{sielie.13}
$$
where we have introduced the quantities
$$
\Theta(r,s,k) = \!\!\!\!\sum_{j_1+\ldots+j_k=s \atop j_1,\ldots,j_k\geq r} 
 \!\!\!\!\frac{(j_1\!+\ldots+\!j_k\!+\!r\!+\!2)
                (j_1\!+\ldots\!+j_{k-1}\!+\!r+\!2)\cdots
                 (j_{1}\!+\!r\!+\!2)}{(j_1+\ldots+j_k)
                                      (j_1+\ldots+j_{k-1})\cdots j_1}
\formula{sielie.13a}
$$
The inequality~\frmref{sielie.13} follows by repeatedly applying
lemma~\lemref{sielie.105} to the non recursive expression for
$E^{(r-1)}_{s-r}$ given in lemma~\lemref{sielie.6}.   We describe this
process in detail.  First estimate
$$
\bignorma{\lie{W^{(r-1)}_{j_1}}\Rmat X_r} 
 \le (j_1+r+2) T_{r-1,j_1} \frac{C_{r-1}^{j_1-1} A}{j_1} 
  \cdot \frac{T_{r-1,r}}{\alpha_r}\cdot\frac{(1+\alpha_r)C_{r-1}^{r-1}
  A}{r}
$$
which follows by the induction hypothesis~\frmref{sielie.71},
from~\frmref{sielie.26} and from lemma~\lemref{sielie.105}.  Remark
that the factor $j_1+r+2$ is the degree of the vector field
$\lie{W^{(r-1)}_{j_1}}\Rmat X_r$.  The same estimate is repeated for
$\lie{W^{(r-1)}_{j_2}}\lie{W^{(r-1)}_{j_1}}\Rmat X_r,\ldots$ thus
yielding
$$
\displaylines{
\quad
\bignorma{\lie{W^{(r-1)}_{j_k}}\cdots\lie{W^{(r-1)}_{j_1}}\Rmat X_r}
 \le \frac{(j_1+\ldots+j_k+r+2)\cdots (j_1+r+2)}{j_k\cdots j_1} \times
\hfill\cr\hfill
  T_{r-1,j_1}\cdots T_{r-1,j_k}\, \frac{T_{r-1,r}}{\alpha_r}
   \cdot \frac{(1\!+\!\alpha_r)C_{r-1}^{s-k-1}A^{k+1}}{r}
\qquad\cr
}
$$
In view of the expression~\frmref{sielie.2} of $E^{(r-1)}_{s-r}$ we
get~\frmref{sielie.13} and~\frmref{sielie.13a}.  Here we need an
estimate for $\Theta(r,s,k)$, that we defer to appendix~\sbsref{B.2}.
Putting $m=2$ in~\frmref{sielie.21} we get
$$
\Theta(r,s-r,k) \lt  2^kr^{k-1} \Bigl(1 +\frac{1}{r}\Bigr)^{k}
  \cbin{\frac{s}{r}-2}{k-1}\ .
\formula{sielie.65}
$$
Furthermore, we isolate the contribution of the quantities
$T_{r-1,\cdot}$ that control the small divisors.  With an appropriate
use of property~(ii) of lemma~\lemref{nrmlie.27} we have
$$
\vcenter{\openup1\jot\halign{
\hfil$\displaystyle{#}$
&$\displaystyle{#}$\hfil
\cr
T_{r-1,j_1}\cdots T_{r-1,j_k}\frac{T_{r-1,r}}{\alpha_r} 
& = \frac{1}{r^2} \left(\frac{\alpha_r}{r^2}\right)^{k-1}
     \frac{r^2 T_{r-1,j_1}}{\alpha_r} \cdots
      \frac{r^2 T_{r-1,j_{k-1}}}{\alpha_r} 
       T_{r-1,j_k} \frac{r^2 T_{r-1,r}}{\alpha_r}
\cr
& \le \frac{1}{r^2} \left(\frac{\alpha_r}{r^2}\right)^{k-1}
   T_{r,s}\ .
\cr
}}
\formula{sielie.113}
$$
Thus, using also $\alpha_r\le 1$ which is true in view of the
definition of $\alpha_r$, we replace~\frmref{sielie.13} with
$$
\eqalign{
\bignorma{E_{s-r}^{(r-1)} \Rmat X_r}
&\le T_{r,s} C_{r-1}^{s-1} \frac{(1+\alpha_r)A}{\alpha_r}
  \sum_{k=1}^{\lfloor s/r\rfloor-1}  
   \Theta(r,s-r,k) \left(\frac{\alpha_r A}{r^2 C_{r-1}}\right)^k
\cr
&\le T_{r,s} C_{r-1}^{s-1} A\cdot 
  \frac{4\left(1+\frac{1}{r}\right)A}{r^2 C_{r-1}}
   \sum_{l=0}^{\lfloor s/r\rfloor-2} \cbin{\frac{s}{r}-2}{l} 
    \left(\frac{2\left(1+\frac{1}{r}\right)A}{r C_{r-1}}\right)^l
\cr
&\lt T_{r,s} C_{r-1}^{s-1} A\cdot 
  \frac{4\left(1+\frac{1}{r}\right)A}{r^2 C_{r-1}}
   \left(1+\frac{2\left(1+\frac{1}{r}\right)A}{r C_{r-1}}\right)^{\frac{s}{r}-2}
    \ .
\cr
}
$$
For $r=1$ the latter formula readily gives~\frmref{sielie.66}.  For
$r\gt 1$ we use $1+\frac{1}{r}\lt 2$ and $\frac{8A}{C_{r-1}}\lt 1$,
which follows from the choice of the sequence $C_r$ in the statement
of the lemma, so that~\frmref{sielie.62} is easily recovered.

\noindent
We come now to the estimate of $\ordnorma{V^{(r)}_s}$.  Here it is
convenient to separate the case $r=1$.  Putting~\frmref{sielie.66} and
the induction hypothesis ~\frmref{sielie.71}
in~\frmref{sielie.8} we get 
$$
\eqalign{
\bignorma{V^{(1)}_s} 
&\lt T_{0,s}\frac{C_{0}^{s-1} A}{s} 
  + T_{1,s} \frac{C_0^{s-1} A}{s} \cdot\frac{8A}{C_{0}}
   \left(1+\frac{4A}{C_{0}}\right)^{s-2}
\cr
&\lt T_{1,s}  \frac{C_0^{s-1} A}{s} 
  \left[1+ \frac{8A}{C_0} \left(1+\frac{4A}{C_{0}}\right)^{s-2}\right]  
\cr
&\lt T_{1,s}  \frac{C_0^{s-1} A}{s} 
  \left(1+\frac{8A}{C_{0}}\right)^{s-1}\ .
\cr
}
$$
Thus we may write
$$
\bignorma{V^{(1)}_s} \lt T_{1,s}  \frac{\hat C_1^{s-1} A}{s}\ ,\quad
 \hat C_1 = 
  C_0 + 8A   \ .
\formula{sielie.64}
$$
For $r\gt 1$ we put~\frmref{sielie.62} in~\frmref{sielie.8}, and we
get
$$
\eqalign{
\bignorma{V^{(r)}_s} 
& \lt T_{r-1,s} \frac{C_{r-1}^{s-1} A}{s} 
   + T_{r,s} \frac{C_{r-1}^{s-1} A}{r s}
  \left(1+\frac{1}{r}\right)^{\frac{s}{r}-2}
\cr
&\lt T_{r,s} \frac{C_{r-1}^{s-1} A}{s} 
  \left[1+\frac{1}{r}
   \left(1+\frac{1}{r}\right)^{\frac{s}{r}-2} \right]
\lt 
 T_{r,s} \frac{C_{r-1}^{s-1} A}{s}
   \left(1+\frac{1}{r}\right)^{\frac{s-1}{r}}
\cr
}
$$
Thus we conclude
$$
\bignorma{V^{(r)}_s} 
\lt  T_{r,s} \frac{\hat C_{r}^{s-1} A}{s}\ ,\quad
 \hat C_r = \left(1+\frac{1}{r}\right)^{1/r} C_{r-1}\ .
\formula{sielie.15}
$$
Now we look for an estimate of $\ordnorma{W^{(r)}_s}$ as given
by~\frmref{sielie.6}.  Recall that we have $s\gt r$, because
$W^{(r)}_r=0$ by construction.  We use~\frmref{sielie.64}
and~\frmref{sielie.15} together with lemma~\lemref{sielie.105}, and
get
$$
\displaylines{
\quad
\bignorma{W^{(r)}_{s}} \le
 T_{r,s}\frac{\hat C_{r}^{s-1} A}{s} + 
  \frac{1}{s}\sum_{k=1}^{\lfloor s/r\rfloor-1}
   \frac{(s+2)(s-r+2)\cdots(s-kr+r+2)}{k!}\times
\hfill\cr\hfill  
\biggl(\frac{r^2 T_{r-1,r}}{\alpha_r}\biggr)^{\!\! k} T_{r,s-kr}  
 \times \biggl(\frac{C_{r-1}^{r-1}A}{r^2} \biggr)^{\!\! k}\,
    {\hat C_r^{s-kr-1}A}\ .
\quad\cr
}
$$
Here we use again the statement~(ii) of lemma~\lemref{nrmlie.27}.  By
the trivial inequality
$$
s-jr+2\le 4r\Bigl(\frac{s}{r} - 1 -j\Bigr)\quad {\rm for}\quad
 0\le j\lt \lfloor s/r\rfloor-1 
$$ 
and remarking that $\hat C_r\gt 4A$ the latter estimate yields
$$
\eqalign{
\bignorma{W^{(r)}_{s}}
&\le T_{r,s}\frac{\hat C_{r}^{s-1} A}{s} \biggl[
 1 + \sum_{k=1}^{\lfloor s/r\rfloor-1}
      \cbin{\frac{s}{r}-1}{k} \left(\frac{4A}{r^2\hat C_{r}}\right)^k
  \biggr]
\cr
&\lt T_{r,s}\frac{\hat C_{r}^{s-1} A}{s}
 \sum_{k=0}^{\lfloor s/r\rfloor-1}
      \cbin{\frac{s}{r}-1}{k} \left(\frac{1}{r^{2}}\right)^k
\cr
&\lt T_{r,s}\frac{\hat C_{r}^{s-1} A}{s}
  \biggl(1+\frac{1}{r^2}\biggr)^{\frac{s-1}{r}}\ .
\cr
}
$$
We conclude
$$
\bignorma{W^{(r)}_{s}} \lt T_{r,s}\frac{C_{r}^{s-1} A}{s}\ ,\quad
 C_{r} = \biggl(1+\frac{1}{r^2}\biggr)^{1/r}\hat C_r\ .
\formula{sielie.67}
$$
In view of~\frmref{sielie.64} and~\frmref{sielie.15} this proves the
claim of the lemma.\endproof

\section{4}{Proof of the main theorem}
Having established the estimate of the iteration lemma~\lemref{sielie.68} on the
sequence of generating vector fields it is now a standard matter to
complete the proof of theorem~\thrref{mainth}.  Hence this section will
be less detailed with respect to the previous ones.

In view of the iteration lemma  we are given an infinite sequence
$\{X_r\}_{r\ge 1}$ of generating vector fields with $X_r$
homogeneous polynomial of degree $r+1$ satisfying
$$
\bignorma{X_r} \le T_{r-1,r} \frac{C_{r-1}^{r-1}A}{r\alpha_r}\ .
$$
By lemma~\lemref{nrmlie.42} we have
$\frac{T_{r-1,r}}{\alpha_r}\le\gamma^r e^{r\Gamma}$ with $\Gamma$ as
in condition~$\tauvet$ and $\gamma$ a constant independent of
$\lambda$.  Moreover, still by lemma~\lemref{sielie.68}, we have
$C_r\le C_{\infty}$, a positive constant independent of $\lambda$.
Thus we have 
$$
\bignorma{X_r} \le \frac{\eta^{r}e^{r\Gamma}}{r} K
\formula{sielie.123}
$$
with positive constants $\eta$ and $K$ independent of $\lambda$.

We refer now to the analytic setting that we recall in
appendix~\sbsref{A.2}.  Every vector field $X_r$ generates a near the
identity transformation $y=\exp\bigl(\lie{X_r})x$ which transforms the
generating sequence $W^{(r-1)}$ into $W^{(r)}$ according to the
algorithm of section~\sbsref{2.3}.  Thus the near the identity
transformation to normal form is generated by the limit $S_X$ of the
sequence of operators $S^{(r)}_X
=\exp(\lie{X_r})\circ\ldots\circ\exp(\lie{X_1})$.  We apply
proposition~\proref{sielie.120} of appendix~\sbsref{A.2}.  In view
of~\frmref{sielie.123}, in a polydisk $\Delta_\rho$ of radius $\rho$
centered at the origin of $\complessi^n$ we have
$$
\bigl|X_r\bigr|_{\rho} 
 \le \bignorma{X_r} \rho^{r+1} 
  \le \frac{\eta^{r}e^{r\Gamma}}{r} K \rho^{r+1}\ ,
$$
where $\bigl|X_r\bigr|_{\rho}$ is the supremum norm.  Thus
condition~\frmref{sielie.121} of proposition~\proref{sielie.120} reads 
$$
\rho \sum_{r\ge 1} \frac{\eta^{r}e^{r\Gamma}}{r}
K \rho^{r} \lt \frac{\rho}{4eK}\ ,
$$
which is true if we take, e.g.,
$\rho\lt\overline\rho=\frac{3}{2}B^{-1}e^{-\Gamma}$ with a constant
$B$ independent of $\lambda$.  Thus proposition~\proref{sielie.120}
applies with, e.g., $\delta=\overline\rho/3$, and we conclude that the
near the identity transformation that gives the map the wanted linear
normal form is analytic at least in a polydisk of radius
$B^{-1}e^{-\Gamma}$, as claimed.  This concludes the proof of the main
theorem.

\vskip 18pt
\centerline{\parfnt Appendix}

\appendix{A}{Technicalities}
This appendix is mainly devoted to recalling some definitions and
properties related to Lie series and Lie transforms that are used in
the paper.  The estimate of the sequence~\frmref{sielie.13a} is
also included.

\subsection{A.1}{Lie series and Lie transforms}
We recall here some formal properties of the Lie series and Lie
transforms as defined by~\frmref{sielie.100} and~\frmref{trslie.1}.
The Lie derivative with respect to the vector field $X$ is defined as
$$
\lie{X}\cdot = \der{}{t}\bigl(\flusso_X^t \cdot\bigr) \Big|_{t=0}
$$
where $\flusso_X^t$ is the flow generated by $X$.  We also recall that
the explicit form of the Lie derivatives for a function $f(x)$ and a
vector field $v(x)$ are, respectively,
$$
\lie{X} f=\sum_{j=1}^{n} X_j\derpar{f}{x_j}
 \quad{\rm and}\quad
 \bigl(\lie{X}v\bigr)_j = \sum_{l=1}^{n} \left(X_l\derpar{v_j}{x_l}
     - v_l\derpar{X_j}{x_l}\right)\ ,
$$
where the l.h.s.~is the $j$-th component of the vector field
$\lie{X}v$.

A non-recursive formula for the operator $E^{X}_s$ defined
by~\frmref{trslie.2} is given by
\lemma{sielie.6}{Let $X=\{X_l\}_{l\gt 1}$ be a sequence of vector
fields and $T_X=\sum_{s\ge 0} E^{X}_s$ as in~\frmref{trslie.2}.  Then for
$s\gt 0$ we have
$$
E^{X}_s = \sum_{k=1}^{s} \sum_{{j_1+\ldots+j_k=s}\atop{j_1,\ldots,j_k\gt 0}}
 \frac{j_k\cdot\ldots\cdot
 j_1}{(j_1+\ldots+j_k)\cdot(j_1+\ldots+j_{k-1})\cdot\ldots\cdot j_1}
  \lie{X_{j_k}}\circ\ldots\circ\lie{X_{j_1}}\ .
\formula{sielie.2}
$$
If moreover $X_1=\ldots=X_{r-1}=0$ then we get the same formula with
the limitation $j_1,\ldots,j_k\ge r$ on the indices $j$.}\endclaim

\proof
The formula may be written down directly by combinatorial
considerations.  A proof by induction goes as follows. We write the
denominator in the equivalent form $ {s(s-j_k)\cdot\ldots\cdot(s-j_k-\ldots-j_2)}\,$. For
$s=1$ both the formula and the recursive definition~\frmref{trslie.2}
give $E^{X}_1=\lie{X_1}$.  Assuming it is true up to $s-1$, for
$l=1,\ldots,s-1$ write the corresponding term in the
r.h.s.~of~\frmref{trslie.2} as
$$
\displaylines{
\frac{l}{s}\lie{X_l}E^{(X)}_{s-l} =
\hfill\cr\hfill 
 \sum_{k=1}^{s-l} \sum_{{j_1+\ldots+j_k=s-l}\atop{j_1,\ldots,j_k\gt 0}}
  \frac{l}{s}\lie{X_l}
   \frac{j_k\cdot\ldots\cdot j_1}{(s-l)(s-l-j_k)\cdot\ldots\cdot(s-l-j_k-\ldots-j_2)}
    \lie{X_{j_k}}\circ\ldots\circ\lie{X_{j_1}}
\cr
}
$$
while for $l=s$ we have the sole term  $\lie{X_s}$.
Renaming the index $l$ as $j_{k+1}$ the corresponding term
$\frac{l}{s}\lie{X_l}$ is included in the sum over the indexes
$j_1,\ldots,j_{k+1}$ with the condition $j_1+\ldots+j_{k+1}=s$, and
writing $k$ in place of $k+1$ the limits of the first sum run from $2$
to $s$, thus changing the sum into
$\sum_{k=2}^{s}\sum_{j_1,\ldots,j_k}\cdots$ so that only the term
$k=1$ is missing.  Adding it means that the second sum contains the
sole term $j=s$, corresponding to  $\lie{X_s}$, which is actually the
missing term in order to recover~\frmref{trslie.2}.\endproof

The Lie series and Lie transform are linear operators acting on the
space of holomorphic functions and of holomorphic vector fields.  They
preserve products between functions and commutators between vector
fields, i.e., if $f,\,g$ are functions and $v,\,w$ are vector fields
then one has
$$
T_X (fg)
=T_X f\cdot T_X g\ ,\quad
T_X \{v,w\} = \bigl\{T_X v, T_X w\bigr\}\ .
\formula{lser.6}
$$
Here, replacing $T_X$ with $\exp\bigl(\lie{X}\bigr)$ gives the corresponding
property for Lie series.  Moreover both operators are invertible.  The
inverse of $\exp\bigl(\lie{X}\bigr)$ is $\exp\bigl(\lie{-X}\bigr)$,
which is a natural fact if one recalls the origin of Lie series as a
solution of an autonomous system of differential equations.

We come now to a remarkable property which justifies the usefulness of
Lie methods in perturbation theory.  We adopt the name {\corsivo
exchange theorem} introduced by Gr\"obner.  Let $f$ be a function and
$v$ be a vector field. {\corsivo Consider the 
transformation $y = T_X x$,
i.e., in coordinates,
$$
y_j = x_j + X_{1,j}(x) 
 +  \left[\frac{1}{2} \lie{X_1} X_{1,j}(x) + X_{2,j}(x)\right] 
  +\ldots
\ ,\quad j=1,\ldots,n\ .
$$
and denote by $\Jmat$ its differential, namely, in
coordinates, the Jacobian matrix with elements $J_{j,k}
=\derpar{y_{j}}{x_k}$.  Then one has
$$
f(y)\big|_{y=T_X x} = 
 \bigl(T_X f\bigr) (x)\ ,\quad
\Jmat^{-1} v(y)\Big|_{y= T_X x} = 
 \bigl(T_X v\bigr)(x)\ .
\formula{trslie.4}
$$}

\noindent
The same statement holds true for a transformation generated by a Lie
series.

The main difference between Lie series and Lie transform is that a
near the identity transformation of the form~\frmref{lser.7} can be
represented via a Lie transform with a suitable sequence of vector
fields.  The same does not hold true for a Lie series, but the same
result may be achieved by using the composition of Lie series.  Again,
see~\dbiref{Giorgilli-2013} for details.

\subsection{A.2}{Analyticity of the compositions of Lie series}
We include here some quantitative estimates which lead to the
convergence of a near the identity coordinate transformation defined by
composition of Lie series.  The results here are basically an
extension of Chapter~1 of Gr\"obner's book~\dbiref{Groebner-1960}.
The reader may also find more details in~\dbiref{Giorgilli-2003.1}
where the Hamiltonian case is dealt with.

We consider a polydisk of radius $\rho$ centered at the origin
$$
\Delta_{\rho} = \bigl\{x\in\complessi^n\>:\> |x|\leq\rho\bigr\}\ .
$$
For a function $f$ and a vector field $X$ it is convenient to
consider the usual supremum norm
$$
|f|_\rho=\sup_{x\in\Delta_{\rho}} |f(x)|
\quad\hbox{and}\quad
|X|_\rho=\sum_j |X_j|_{\rho}\ .
$$

The following lemma states the Cauchy inequality for Lie derivatives
in $\Delta_{\rho}$.

\lemma{norms:liederivative}{For $0\leq \delta' <\delta<\rho$ the following
inequalities apply to Lie derivatives:
$$
|\lie{X}{f}|_{\rho-\delta} \leq \frac{1}{\delta-\delta'} |X|_\rho
 |f|_{\rho-\delta'}\ ,
$$
and, for $s\geq1$,
$$
|\lie{X}^{s}{f}|_{\rho-\delta} \leq \frac{s!}{e} \left(\frac{e
 |X|_\rho}{\delta}\right)^s |f|_{\rho}\ .
$$
}\endclaim

\noindent
The proof is a straightforward adaptation of Cauchy estimates for
derivatives of analytic functions.  The reader will be able to
reproduce it himself, perhaps with the help of the proof for the
Hamiltonian case in~\dbiref{Giorgilli-2003.1}.

The following proposition essentially restates Cauchy's existence
theorem of solution of differential equations in the analytic case.

\proposition{norms:explie}{Let $X$ be an analytic vector field on the
interior of $\Delta_{\rho}$ and bounded on $\Delta_{\rho}$.  Then, for
every positive $\delta\leq\rho/2$ the following statement holds true:
if $|X|_{\rho}\lt\delta(e-1)/e^2$ then the map
$$
\phi(x) = \exp(\lie{X})x
\formula{trasf.1}
$$
is bianalytic and satisfies
$$
|\exp(\lie{X})x_j - x_j|_{\rho-\delta}\leq\frac{\delta}{e^2}\ ,
\formula{stime:trasf}
$$
so that we have
$$
\Delta_{\rho-2\delta} \subset \phi(\Delta_{\rho-\delta}) \subset \Delta_{\rho}
\ ,\quad
\Delta_{\rho-2\delta} \subset \phi^{-1}(\Delta_{\rho-\delta}) \subset \Delta_{\rho}\ .
$$}\endclaim
\proof{Using lemma~\lemref{norms:liederivative}, compute
$$
\displaylines{
\qquad
|\exp(\lie{X})x_j - x_j|_{\rho-\delta}
 = \sum_{s\geq1} \Bigl| \frac{\lie{X}^{s-1}}{s!} X_j\Bigr|_{\rho-\delta}
\hfill\cr\hfill
 \leq \sum_{s\geq1} \frac{1}{es} \Bigl( \frac{e
  |X|_\rho}{\delta}\Bigr)^{s-1} |X|_\rho
  \leq -\frac{\delta}{e^2} \ln\left(1-\frac{e|X|_\rho}{\delta}\right)
   \leq\frac{\delta}{e^2}\ ,
\qquad\cr
}
$$
namely~\frmref{stime:trasf}.  The same estimates apply to the inverse
$\exp\bigl(\lie{-X}\bigr)$.}\endproof

Let now $\{X_1,\,X_2,\ldots\}$ be a sequence of analytic vector
fields.  Consider the sequence of transformations
$\{S^{(0)}_X,\,S^{(1)}_X,\,S^{(2)}_X,\ldots\}$ recursively defined as
$$
S^{(0)}_X = \uno\ ,\quad
S^{(r)}_X = \exp\bigl(\lie{X_r}\bigr) \circ S^{(r-1)}_X\ .
\formula{lser.8r}
$$
with inverse
$$
\tilde S^{(0)}_X = \uno\ ,\quad
\tilde S^{(r)}_X = \tilde S^{(r-1)}_X\circ \exp\bigl(\lie{-X_r}\bigr) \ .
\formula{lser.8rinv}
$$
We may well consider in formal sense the composition of Lie series as
the limit
$$
S_X = \ldots\circ\exp\bigl(\lie{X_r}\bigr)\circ\ldots\circ
     \exp\bigl(\lie{X_2}\bigr)\circ
      \exp\bigl(\lie{X_1}\bigr)\ ,
\formula{lser.8}
$$
with a similar formal limit for the sequence $\tilde S_X^{(r)}$ giving
the inverse $S^{-1}_X$.

\proposition{sielie.120}{Let the sequence $X=\{X_r\}_{r\ge 1}$ of vector
fields be analytic in $\Delta_{\rho}\,$. Then for every positive
$\delta\lt\rho/2$ the following statement holds true: if the convergence
condition
$$
\sum_{r\ge 1} \bigl|X_r\bigr|_{\rho} \lt \frac{\rho}{4e}
\formula{sielie.121}
$$
is satisfied then the sequence
$\{S^{(0)}_X,\,S^{(1)}_X,\,S^{(2)}_X,\ldots\}$ of maps defined
by~\frmref{lser.8r} is bianalytic and converges to a bianalytic
map $S_X$ satisfying
$$
\Delta_{\rho-2\delta} \subset S_X(\Delta_{\rho-\delta}) \subset \Delta_{\rho}
\ ,\quad
\Delta_{\rho-2\delta} \subset
S_X^{-1}(\Delta_{\rho-\delta}) \subset \Delta_{\rho}\ .
\formula{sielie.122}
$$}\endclaim

\proof
Consider the sequence of positive numbers
$$
\delta_s = 
 \frac{|X_s|_{\rho}}{\sum_{r\ge 1}|X_r|_{\rho}} \delta\ ,\quad
  s\ge 1
$$
so that $\sum_{s\ge 1}\delta_s=\delta$.  By definition we have
$|X_s|_{\rho}\lt \frac{\delta_s}{2e}$, so that
proposition~\proref{norms:explie} applies to every vector field
$X_s\,$.  Consider also the two sequences
$$
\rho^-_0 = \rho^+_0 = \rho-\delta\ ,\quad
\rho^-_s = \rho^-_{s-1}-\delta_s\ ,\quad
\rho^+_s = \rho^+_{s-1}+\delta_s\ ,\quad s\ge 1\ ,
$$
so that we have
$\rho^-_s\mathop{\longrightarrow}\limits_{s\to+\infty}\rho-2\delta$
and $\rho^+_s\mathop{\longrightarrow}\limits_{s\to+\infty}\rho$.
Applying proposition~\proref{norms:explie} we see that the map
$S^{(1)}$ is analytic and satisfies $\Delta_{\rho^-_1}\subset
S^{(1)}\Delta_{\rho_0}\subset\Delta_{\rho^+_1}$, and that its inverse
is analytic and satisfies $\Delta_{\rho^-_1}\subset
S^{(1)}\Delta_{\rho_0}\subset\Delta_{\rho^+_1}\,$. Proceeding by
induction it is straightforward to check that the same claim holds
true for $S^{(r)}$ and $\tilde S^{(r)}$ for every $r\gt 1$, with
$\rho^-_1$ and $\rho^+_1$ replaced by $\rho^-_r$ and $\rho^+_r$,
respectively.  Moreover the inequalities
$\bigl|S^{(r)}x-S^{(r-1)x}\bigr|\lt \delta_r$ and $\bigl|\tilde
S^{(r)}x-\tilde S^{(r-1)x}\bigr|\lt \delta_r$ hold true.  Thus the
sequences $S^{(r)}$ and $\tilde S^{(r)}$ of maps is uniformly
convergent in any compact subset of
$\Delta_{\rho_0}=\Delta_{\rho-\delta}$, and by Weierstrass theorem 
the limits $S_x$ and $S^{-1}_X$ do exist and are analytic maps
satisfying~\frmref{sielie.122}.\endproof

\subsection{B.2}{Estimate of the sequence~\frmref{sielie.13a}}
We introduce a parameter
$m\ge 0$ which is arbitrary but fixed and rewrite the
sequence~\frmref{sielie.13a} as
$$
\Theta(r,s,k) = \!\!\!\!\sum_{j_1+\ldots+j_k=s \atop j_1,\ldots,j_k\geq r} 
 \!\!\!\!\frac{(j_1\!+\ldots+\!j_k\!+\!r\!+\!m)
                (j_1\!+\ldots\!+j_{k-1}\!+\!r+\!m)\cdots
                 (j_{1}\!+\!r\!+\!m)}{(j_1+\ldots+j_k)
                                      (j_1+\ldots+j_{k-1})\cdots j_1}
\formula{sielie.12}
$$
with $s\ge kr$ and $1\le k\le \lfloor s/r\rfloor$.  
With a straightforward calculation we get the uniform estimate
$$
\eqalign{
&\frac{(j_1\!+\ldots+\!j_k\!+\!r\!+\!m)
                (j_1\!+\ldots\!+j_{k-1}\!+\!r+\!m)\cdots
                 (j_{1}\!+\!r\!+\!m)}{(j_1+\ldots+j_k)
                                      (j_1+\ldots+j_{k-1})\cdots j_1}
\cr&\qquad
 =  \Bigl(1 +\frac{r+m}{j_1+\ldots+j_k}\Bigr)\cdot
  \Bigl(1 +\frac{r+m}{j_1+\ldots+j_{k-1}}\Bigr)\cdots
    \Bigl(1 +\frac{r+m}{j_1}\Bigr)
\cr&\qquad
\lt \Bigl(2 +\frac{m}{r}\Bigr)^{k} \ ,
\cr
}
$$
which holds true because $j_1,\ldots,j_k\ge r$.  Thus we get
$$
\Theta(r,s,k) = \!\!\!\!\sum_{j_1+\ldots+j_k=s \atop j_1,\ldots,j_k\geq r}
 \Bigl(2 +\frac{m}{r}\Bigr)^{k} \ .
$$
Writing $j_1=r+q_1,\ldots,j_k=r+q_k$ the number of terms in the sum is
equal to  the number of non negative integer vectors
$q=(q_1,\ldots,q_n)$ with $|q|= s-kr$, which is known to be
$$
 \cbin{s-kr+k-1}{k-1}\le  r^{k-1}\cbin{\frac{s}{r}-1}{k-1}\ .
$$
We conclude
$$
\Theta(r,s,k) \lt 
 r^{k-1} \Bigl(2 +\frac{m}{r}\Bigr)^{k}
  \cbin{\frac{s}{r}-1}{k-1}\ .
\formula{sielie.21}
$$

{\bf Acknowledgement}
\par\nobreak\noindent
This research is partially supported by PRIN-MIUR 2010JJ4KPA ``Teorie
geometriche e analitiche dei sistemi Hamiltoniani in dimensioni finite
e infinite''.  M.~S. has also been partially supported by an FSR
Incoming Post-doctoral Fellowship of the Acad\'emie universitaire
Louvain, co-funded by the Marie Curie Actions of the European
Commission.

\references

\bye